\documentclass{article}

\usepackage{amsmath,amsthm,amsfonts,amssymb,latexsym,epsfig}
\usepackage{url}

\theoremstyle{plain}
\newtheorem{theorem}{Theorem}
\newtheorem{corollary}{Corollary}

\begin{document}

\title{Recurrent Proofs of the Irrationality of Certain Trigonometric Values}
\author{Li Zhou and Lubomir Markov}
\date{}
\maketitle

In this note we exploit recurrences of integrals to give new elementary proofs of the irrationality of $\tan r$ for $r\in\mathbb{Q}\setminus\{0\}$ and $\cos r$ for $r^2\in\mathbb{Q}\setminus\{0\}$. We also discuss applications of our technique to simpler irrationality proofs such as those for $\pi$, $\pi^2$, and certain values of exponential and hyperbolic functions.

\section{Irrationality of $\tan r$ for $r\in\mathbb{Q}\setminus\{0\}$.}

For a nonzero rational $r$, the irrationality of $\tan r$ was first proved by J. H. Lambert in 1761 by means of continued fractions \cite[pp.\ 129--146]{bbb}. We now present a new direct proof using a recurrence for an integral.

\begin{theorem}\label{th1} 
$\tan r$ is irrational for nonzero rational $r$.
\end{theorem}

\begin{proof}
The irrationality of $\pi$ will be a by-product of this proof, so we start by supposing that $r\in \mathbb{Q}\setminus \{k\pi : k\in\mathbb{Z}\}$. Write $r=a/b$ with $a,b\in\mathbb{Z}$ and assume that $\tan (r/2)=p/q$ with $p,q\in \mathbb{Z}$. For $n\ge 0$, let $f_n(x)=(rx-x^2)^n/n!$ and $I_n=\int_0^r f_n(x)\sin x \, dx$. Then $b^n I_n\rightarrow\ 0$ as $n\rightarrow\infty$, $I_0=1-\cos r$, and $I_1=2(1-\cos r)-r\sin r$. Integrating by parts twice, we get that for $n\ge 2$, 
\begin{equation}\label{rec0}
I_{n}=- \int_0^{r}f''_n (x) \sin x \,dx= (4n-2)I_{n-1}-{r}^2 I_{n-2}.\end{equation}
By inducting on $n$ using (\ref{rec0}), we see that for $n\ge 0$, $I_n=u_n(1-\cos r)+v_n\sin r$, where $u_n$ and $v_n$ are polynomials in $r$ with integer coefficients and degrees at most $n$. Moreover, if two consecutive terms of the sequence $\langle  I_n\rangle$ are $0$, then (\ref{rec0}) forces all terms of $\langle  I_n\rangle$ to be $0$, and in particular $I_0=0$, a contradiction. Hence $\langle I_n \rangle$ has infinitely many nonzero terms. Therefore, we can pick a large enough $n$ so that $b^n q\csc r I_n=b^n q[u_n\tan(r/2)+v_n]$ is a nonzero integer in $(-1,1)$, a contradiction. 

Notice that $\tan(\pi/4)=1$, so $\pi/2\notin\mathbb{Q}$, which implies that $\mathbb{Q}\setminus \{k\pi : k\in\mathbb{Z}\}=\mathbb{Q}\setminus\{0\}$. Thus we have proved that $\tan(r/2)\notin\mathbb{Q}$ for all $r\in \mathbb{Q}\setminus\{0\}$.
\end{proof}

A closer inspection of our proof reveals that $u_n$ and $v_n/r$ are polynomials in $r^2$. So a slightly stronger conclusion can be squeezed out of the proof, namely that $(\tan r)/r$ is irrational whenever $r^2\in \mathbb{Q}\setminus\{0\}$. This stronger result was first established through a different elementary approach by Inkeri \cite{ink}.
 
\section{Simpler Proofs of Irrationality.}

From the previous proof we see that recurrence is a double-edged sword. It is sharp and swift in showing that a sequence is integer-valued and has an infinite nonzero subsequence. We can also use recurrence to give similar proofs of the irrationality of $\pi$, $\pi^2$, $e^r$, etc.. However, in these easier cases, the corresponding sequences of integrals are positive, so there is no need to argue for the existence of a nonzero subsequence. Consequently the proofs can be really short and charming. For a flavor of it, we present a direct proof of the irrationality of $\pi$ which is even shorter than the celebrated one-page proof given by Niven \cite{n01}.
   
\begin{theorem}\label{th2}
$\pi$ is irrational.
\end{theorem}

\begin{proof}
Assume that $\pi = a/b$ with $a,b\in \mathbb{N}$. Let $f_n(x)=(\pi x-x^2)^n/n!$ and $I_n=\int_0^{\pi} f_n(x)\sin x\, dx$. Then $b^n I_n\rightarrow 0$ as $n\rightarrow\infty$, $I_0=2$, and $I_1=4$. Replacing $r$ by $\pi$ in (\ref{rec0}) we get $I_{n}= (4n-2)I_{n-1}-{\pi}^2 I_{n-2}$ for $n\ge 2$.
By induction on $n$ using this recurrence, we see that for $n\ge 0$, $I_n$ is a polynomial in $\pi$ with integer coefficients and degree at most $n$. Hence for a large enough $n$, $b^nI_n$ is an integer in $(0,1)$, a contradiction.
\end{proof} 

Notice that the terms of $\langle I_n \rangle$ are really polynomials in $\pi^2$, so our proof only needs very minor changes to show the stronger conclusion that $\pi^2$ is irrational. In fact in \cite{sch}, Schr\"oder presented a very similar proof by recurrence of the irrationality of $\pi^2$. 

For the irrationality of $e^r$ for nonzero rational $r$, the interested reader can imitate the process with the sequence $I_n=\int_0^r f_n(x)e^x\, dx$, where $f_n(x)=(rx-x^2)^n/n!$ as well.   

\section{Irrationality of $\cos r$ for $r^2\in\mathbb{Q}\setminus\{0\}$.}

Next we turn our attention to the cosine function. The classical elementary proof of the irrationality of $\cos r$ for nonzero rational $r$ was given by Niven \cite[Theorem 2.5, pp.\ 16--19]{n02}. Niven's proof can be slightly modified to show that $\cos r$ is irrational whenever $r^2 \in\mathbb{Q}\setminus\{0\}$, as observed by Inkeri \cite{ink}.  We now use the full force of recurrence to give a proof which is more direct than Inkeri's modification.
 
\begin{theorem}\label{th3}
If $r^2 \in\mathbb{Q}\setminus\{0\}$ then $\cos r$ is irrational.
\end{theorem}

\begin{proof}
Assume that $r^2=a/b$ with $a,b\in\mathbb{Z}\setminus\{0\}$ and $\cos  r=p/q$ with $p,q\in \mathbb{Z}$. For $n\ge 0$, let $f_{n}(z)=(r^2z^2-z^4)^n/n!$, $I_{n}=\int_0^r f_{n}(z)\sin (r-z)\, dz$, $J_{n}=\int_0^r zf_{n}(z)\cos (r-z) \,dz$,
$K_{n}=\int_0^r z^2f_{n}(z)\sin (r-z)\, dz$, and $L_{n}=\int_0^r z^3f_{n}(z)\cos (r-z)\, dz$. Then $b^{2n+1}$ times each of the four integrals approaches $0$ as $n\rightarrow \infty$. Direct integration yields $I_0=1-\cos r=J_0$, $K_0=r^2-2+2\cos r$, and $L_0=3K_0$. Integrating each integral by parts once, we get that for $n\ge 1$,
\begin{eqnarray}\label{rec1}
I_{n} &=& 4L_{n-1}-2r^2J_{n-1},\\
\label{rec2}
J_{n} &=& (4n+1)I_{n}-2r^2K_{n-1},\\
\label{rec3}
K_{n} &=& -(4n+2)J_{n}+2r^2L_{n-1},\\
\label{rec4}
L_{n} &=& (4n+3)K_{n}+2nr^2I_{n}-2r^4K_{n-1}.
\end{eqnarray}
Induction on $n$ in these recurrences implies that for $n\ge 0$, the four sequences have the form $u_n+v_n\cos r$, where $u_n$ and $v_n$ are polynomials in $r^2$ with integer coefficients and degrees at most $2n+1$. Moreover, suppose that $I_m=J_m=K_m=L_m=0$ for some $m\ge 1$. Then (\ref{rec2}) and (\ref{rec3}) imply that $K_{m-1}=0$ and $L_{m-1}=0$. Thus (\ref{rec1}) yields $J_{m-1}=0$. Eliminating $K_{n-1}$ from (\ref{rec2}) and (\ref{rec4}), we see that for $n\ge 1$, $I_n$ can be expressed in terms of $J_n$, $K_n$, and $L_n$; and keep in mind also that $I_0=J_0$. Hence $I_{m-1}=0$. By this argument of infinite descent we conclude that $I_0=J_0=K_0=L_0=0$, which contradicts the fact that $2I_{0}+K_0=r^2\ne 0$. Therefore, at least one of the four sequences has infinitely many nonzero terms. Pick such a sequence and a large enough $n$ so that $b^{2n+1} q$ times the corresponding integral is a nonzero number in $(-1,1)$, while it also has the form $b^{2n+1}q(u_n+v_n\cos r)$ which is an integer, a contradiction. Thus $\cos r\notin\mathbb{Q}$ whenever $r^2\in\mathbb{Q}\setminus\{0\}$.
\end{proof}

\begin{corollary}
${\pi}^2$ is irrational. Also, if $r^2\in\mathbb{Q}\setminus\{0\}$ then $\sin^2 r$, $\cos^2 r$, and $\tan^2 r$ are all irrational.
\end{corollary}

\begin{proof}
The claims follow immediately from Theorem \ref{th3} and the identities $\cos \pi =-1$ and $\cos 2r =1-2\sin ^2 r=2\cos^2 r-1=(1-\tan^2 r)/(1+\tan ^2 r)$.
\end{proof}

Finally, the observant reader perhaps has noticed that our proof of Theorem \ref{th3} allows $r^2$ to be a negative rational number. Since $\cosh r=\cos(ir)$, all the analogous statements about hyperbolic functions are included in our results. The skeptical reader is invited to work out the details by substituting $r=it$ and $z=iy$ into $I_n$, $J_n$, $K_n$, and $L_n$. In fact, the resulting real integrals are nonzero and thus the proof is shorter, because the argument of infinite descent is not needed.

\paragraph{Acknowledgments.} We are grateful to our friend Prof. Paul Yiu whose comments on an earlier draft led to a substantial improvement of the paper. We would also like to thank the anonymous referees for their suggestions and for informing us of the references \cite{ink} and \cite{sch}.

\bigskip

\noindent\textit{Department of Mathematics, Polk State College,
Winter Haven, FL 33881\\
lzhou@polk.edu}

\bigskip

\noindent\textit{Department of Mathematics and CS, Barry University, Miami Shores, FL 33161\\
lmarkov@mail.barry.edu}


\begin{thebibliography}{9}
\bibitem{bbb} L. Berggren, J. Borwein, and P. Borwein, eds., \textit{Pi: A Source Book}, 3rd ed., Springer, New York, 2004.
\bibitem{ink} K. Inkeri, Er\"aiden transsendenttifunktioiden arvojen irrationaalisuudesta, \emph{Arkhimedes}, no. 1 (1965) 15--21.
\bibitem{n01} I. Niven, A simple proof that $\pi$ is irrational, \emph{Bull. Amer. Math. Soc.} \textbf{53} (1947) 509.
\bibitem{n02} ---------, \textit{Irrational Numbers}, Mathematical Association of America, Washington, DC, 1956.
\bibitem{sch} E. M. Schr\"oder, Zur irrationalit\"at von $\pi^2$ und $\pi$, \emph{Mitt. Math. Ges. Hamburg} \textbf{13} (1993) 249.
\end{thebibliography}
\end{document}